\documentclass[twocolumn]{autart}

\usepackage{tikz}
\usepackage{latexsym, amssymb}
\usepackage{graphicx}
\usepackage{amsmath, amsbsy}
\usepackage{amsopn, amstext}
\usepackage{hyperref}
\usepackage{appendix}
\usepackage{cancel, color}
\usepackage{epstopdf}
\usepackage[noend]{algpseudocode}
\usepackage{algorithmicx,algorithm}

\newtheorem{dfn}[thm]{Definition}
\newtheorem{prp}[thm]{Proposition}
\newtheorem{exa}[thm]{Example}

\begin{document}

\begin{frontmatter}

\title{A Strategic Learning Algorithm for State-based Games \thanksref{footnoteinfo}}

\thanks[footnoteinfo]{This work is supported partly by the National Natural Science Foundation of China (NSFC) under Grants 61473099, 61773371, 61733018 and 61333001. Corresponding author: Fenghua He. Tel.: +86 0451-86402947; fax.: +86 0451-86414580.}

\author{Changxi Li}\dag\ead{changxi1989@163.com},\author{Yu Xing}\ddag\ead{yxing@amss.ac.cn}, \author{Fenghua He}\dag\ead{hefenghua@gmail.com}, and \author{Daizhan Cheng}\ddag\ead{dcheng@iss.ac.cn}

\address{\dag Control and Simulation Center, Harbin Institute of Technology, Harbin 150001, P.~R.~China\\
         \ddag Academy of Mathematics and Systems Science, Chinese Academy of Sciences, Beijing 100190, P.R.China}

\begin{abstract}
Learning algorithm design  for  \emph{state-based games} is investigated.  A heuristic uncoupled learning algorithm, which  is  a two memory  better reply with inertia dynamics, is proposed.  Under certain reasonable conditions it is proved that for any initial state, if all agents in the state-based game  follow the proposed  learning algorithm,  the action state pair converges almost surely to an action invariant set of recurrent state equilibria. The design relies on  global and local searches with finite memory, inertia, and randomness. Finally,  existence of time-efficient universal learning algorithm is studied. A class of state-based games is presented to show that there is  no universal learning algorithm converging  to a recurrent state equilibrium.
\end{abstract}


\end{frontmatter}

\section{Introduction}
Many systems, such as biological  networks, social networks \cite{jr56}, and engineering systems \cite{fb18}, can be described as a collection of interacting subsystems, which causes local decisions using local information \cite{jr15}.
To ensure the emergence of desirable collective behavior by  designing proper local control strategies is the core mission in such systems.
Game-theoretical method is becoming an appealing tool in control of the above  systems as it provides a modularized design architecture, i.e.  the interaction structure and learning algorithms can be designed separately \cite{jr15},\cite{co17}. Some outstanding works include: (i) consensus/synchronization of multi-agent systems \cite{jr15}; (ii) distributed optimization \cite{by10}; (iii) optimization in energy \cite{wss12} and transportation networks \cite{wan13}, just to name a few.



State-based games, an extended model in game-theoretic control, were  proposed in \cite{jm12}. In fact, the idea of state-based games can be traced back to \cite{hp04} (Section $9$, Conclusion). Since then state-based games have shown their strong vitality in many fields, such as  achieving  Pareto optimality \cite{jh14},  realizing  cooperative coverage  in unknown environment \cite{srw14}, and  solving distributed economic problem in smart grid \cite{yl16}. Particularly, a completely uncoupled learning algorithm for general  games is designed for the first time using the theory of state-based games and regular perturbed Markov chain \cite{hp09}.

%
Compared with traditional game-theoretical framework, state-based games provide an additional degree of freedom, which is called \emph{state},  to help coordinate group behavior. The underlying ``state" has a variety of interpretations ranging from a dummy agent \cite{jm12} or external environment \cite{hp04} to  real agents with unknown dynamics or dynamics for equilibrium selection \cite{bsr12},\cite{jr17}. Since the additional degree of freedom is provided to help coordinate group behavior, state-based game is a  useful extended model in game-theoretic control.

One of the core challenges in applying state-based game method to game-theoretic control is to  design a strategic learning algorithm  which can converge to the equilibria of  state-based games. Although \cite{jm12} proposed a finite memory learning algorithm for state-based potential games, to our best knowledge, there is no strategic learning algorithm for general state-based games.
The purpose of this paper is to design a heuristic algorithm for general  state-based games.


The main contribution of this paper is the designed  two memory strategic learning  algorithm for general  state-based games. The designed algorithm relies on global and local searches using two memory information, inertia, and randomness. Under certain reasonable conditions it is proved that the algorithm  converges almost surely to a recurrent state equilibrium of state-based games,  which is a generalized Nash equilibrium. Finally, to investigate the existence of universal learning algorithm, a class of state-based games is presented, and for such state-based games there is no universal learning algorithm converging  to a recurrent state equilibrium.

%

The rest of this paper is organized as follows: Section $2$ provides some preliminaries, including the formal definition of state-based games, recurrent state equilibrium, state-based potential games, and the theory of learning in state-based games. Section $3$ focuses  on the design of a learning algorithm for general state-based games. Section $4$ considers the existence of a universal learning algorithm. A brief conclusion is given in Section $5$. Appendix contains three parts. First part reveals the Markov chain induced by the proposed learning algorithm. Some lemmas used in the proof of the convergence of the proposed learning algorithm are provided in  Second part. The convergence of the proposed learning algorithm  is proved in  last part.

\section{Preliminaries}
\subsection{State-based games}
\begin{dfn}\cite{jm12} \emph{(State-based game)}
A finite state-based game is a quintuple $\mathcal{G}=\left\{N,\{A_i\},\{c_i\},X,P\right\}$, where
\begin{enumerate}
\item $N=\{1,2,\cdots,n\}$ is the set of agents;
\item $A_i=\{1,2,\cdots,k_i\}$ is the set of actions of agent $i$;
\item $c_i:~A\times X\rightarrow \mathbb{R}$ is the payoff function of agent $i\in N$, where $A=\prod_{i=1}^nA_i$ is the action profile set, and  $\prod$ is the Cartesian product;
\item $X=\{1,2,\cdots,m\}$ is the set of underlying finite state;
\item $P:A\times X\rightarrow \Delta(X)$ is the Markovian state transition function, where $\Delta(X)$ denotes the set of probability distributions over  the finite state space $X$.
\end{enumerate}
\end{dfn}

\begin{figure}[!hbtp]
\centering
\includegraphics[ width = 3.3 in]{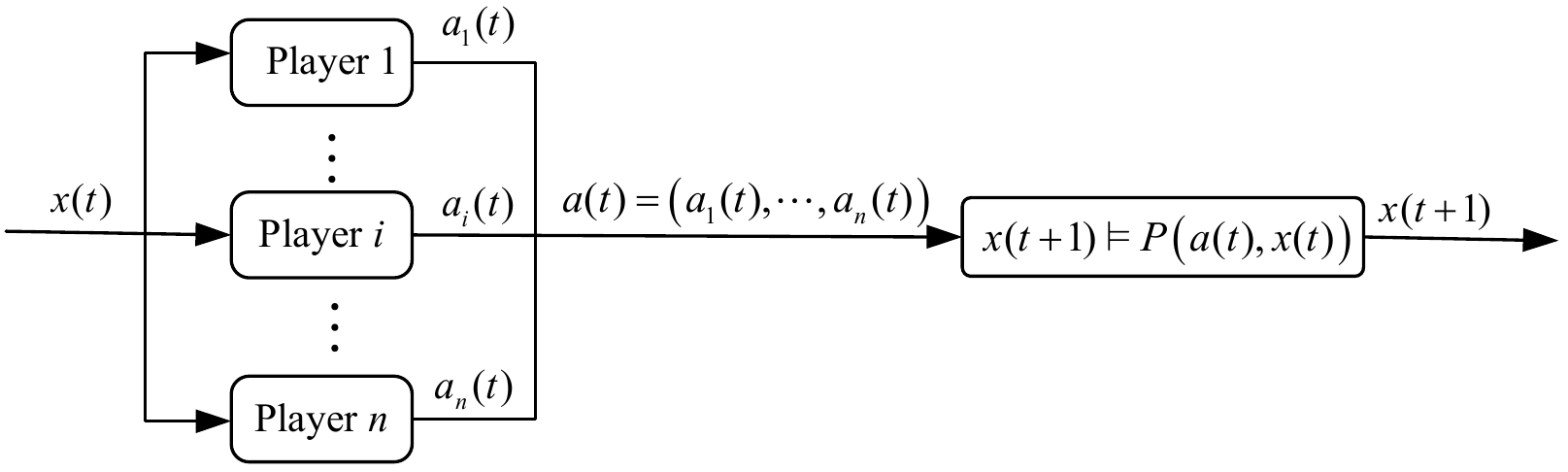}
\caption{Dynamics of State-based Games}
\label{fig00.1}
\end{figure}

When a state-based game is played repeatedly, a sequence of states
$$x(0), x(1),\cdots,x(t),\cdots$$
and a sequence of joint actions
$$a(0), a(1),\cdots,a(t),\cdots$$
are generated. $[a(t),x(t)]\in A\times X$ is referred to the action state pair at time $t$. We give a rough description on how the action state pair evolves.  The sequence of  action profiles is produced using some specified decision algorithm. Suppose the current state is $x(t)$, and the action taken by all agent at time $t$ is $a(t)$, then $x(t+1)$ is generated  by the state transition function $P(a(t),x(t))$, i. e., the ensuing state is selected  randomly according  to the probability distribution $P(a(t),x(t))$. The dynamics of state-based games can be described as in Fig. \ref{fig00.1}, where `$\vDash$' signifies that the ensuing state $x(k+1)$ is selected according to the probability distribution $P(a(k),x(k))$.

Denote by $X(a|x)\subseteq X$ the set of reachable states starting from initial state $x$ driven by an invariant action $a$. That is to say, a state $y\in X(a|x)$ if and only if there exists a time $t_{y}>0$ such that
$$\textbf{Pr}[x(t_{y})=y]>0,$$
conditioned on the events $x(0)=x$ and $x(k+1)\vDash P(a,x(k))$ for all $k\in \{0,1,\cdots,t_y-1\}$. The transition process can be illustrated as
$$x\stackrel{a}{\longrightarrow}x(1)\stackrel{a}{\longrightarrow}~\stackrel{a}{\cdots}~\stackrel{a}{\longrightarrow}x(t_y-1)\stackrel{a}{\longrightarrow}x(t_y)=y.$$
\begin{rem}
As pointed in \cite{jm12} (Section 3.4), the model of state-based games is a simplification of Markov games \cite{shp53}. In state-based games each agent is myopic (seeks to optimize the current payoffs), while in markov games every agent seeks to optimize a discounted sum of future payoffs.
\end{rem}

As a generalization of Nash equilibrium, the equilibrium in state-based games is called the \emph{recurrent state equilibrium} (RSE).
\begin{dfn} \cite{jm12} \emph{(Recurrent state equilibrium)}
Consider a state-based game $\mathcal{G}=\left\{N,\{A_i\}_{i\in N},\{c_i\}_{i\in N},X,P\right\}.$
The action state pair $[a^*,x^*]$ is a recurrent state equilibrium with respect to the state
transition process $P(\cdot)$ if the following two conditions are satisfied:
\begin{enumerate}
\item The state $x^*$ satisfies $x^*\in X(a^*|x)$ for every state $x\in X(a^*|x^\ast)$;
\item For each agent $i\in N$ and every state $x\in X(a^*|x^*)$,
$$c_i(a_i^*,a_{-i}^*,x)\geqslant c_i(a_i,a_{-i}^*,x),~\forall a_i\in A_i.$$
\end{enumerate}
\end{dfn}
Denote  $P(a;\cdot,\cdot)$ the probability transition matrix of a joint action $a\in A$ in a state-based game  $\mathcal{G}$.
The first condition means that if the action state pair $[a^*,x^*]$ is a recurrent state equilibrium, then  $X(a^*|x^*)$ is a recurrent class of the Markov chain $P(a^*;\cdot,\cdot)$ starting from the initial state $x^*$. The second condition implies that $a^*$ is a pure Nash equilibrium of state invariant game $G_x=\{N,A_i,c_i(\cdot,x)\}$ for every state $x\in X(a^*|x^*).$

Consider two action state pairs $[a,x]$ and $[b,y]$. $[a,x]$ and $[b,y]$ are called \emph{equivalent} if the following three conditions are satisfied:
\emph{i)} $a=b$, \emph{ii)} $[a,x]$ is a recurrent state equilibrium, and \emph{iii)} $y\in X(a|x)$. Use the notation $[a,x]\sim [b,y]$ to represent that $[a,x]$ and $[b,y]$ are equivalent. Otherwise, it is denoted by $[a,x]\nsim [b,y]$. It is easy to verify that $\sim$ is an equivalence relation. Denote  $$R(a,x):=\big\{[a,y]:~[a,y]\sim[a,x]\big\}.$$
We call $R(a,x)$ a recurrent state equilibrium set generated by the recurrent state equilibrium $[a,x]$.
\begin{exa}\label{ex1} Consider the following state-based game with $N=\{1,2\},$ $A_1=A_2=\{1,2\},$ $X=\{1,2,3\}$. The game $G_x$ is a coordination game, prisoner's dilemma game, and matching pennies game when $x=1,2$, and $3$, respectively. The payoff matrices are shown as follows.
\begin{table}[!htbp] 
\centering
\caption{Payoff Bi-Matrix for $x=1$ of Example \ref{ex1} (coordination game)  \label{Tab1.1}}
\begin{tabular}{ccc}
\hline Agent $1\backslash$Agent $2$ &$1$&$2$\\
\hline $1$&$(4,~4)$&$(1,~3)$\\
 $2$&$(3,~1)$&$(2,~2)$\\
\hline
\end{tabular}
\end{table}
\begin{table}[!htbp]
\centering
\caption{Payoff Bi-Matrix for $x=2$ of Example \ref{ex1} (prisoner's dilemma game)  \label{Tab1.2}}
\begin{tabular}{ccc}
\hline Agent $1\backslash$Agent $2$ &$1$&$2$\\
\hline $1$&$(2,~2)$&$(0,~3)$\\
 $2$&$(3,~0)$&$(1,~1)$\\
\hline
\end{tabular}
\end{table}
\begin{table}[!htbp]
\centering
\caption{Payoff Bi-Matrix for $x=3$ of Example \ref{ex1} (matching pennies game) \label{Tab1.3}}
\begin{tabular}{ccc}
\hline Agent $1\backslash$Agent $2$ &$1$&$2$\\
\hline $1$&$(-1,~1)$&$(1,~-1)$\\
 $2$&$(1,~-1)$&$(-1,~1)$\\
\hline
\end{tabular}
\end{table}

The state transition process is shown in Fig.  \ref{fig1.1}.
\begin{figure}[!hbtp]
    \centering
    \includegraphics[height = 6.5cm, width = 6.5cm]{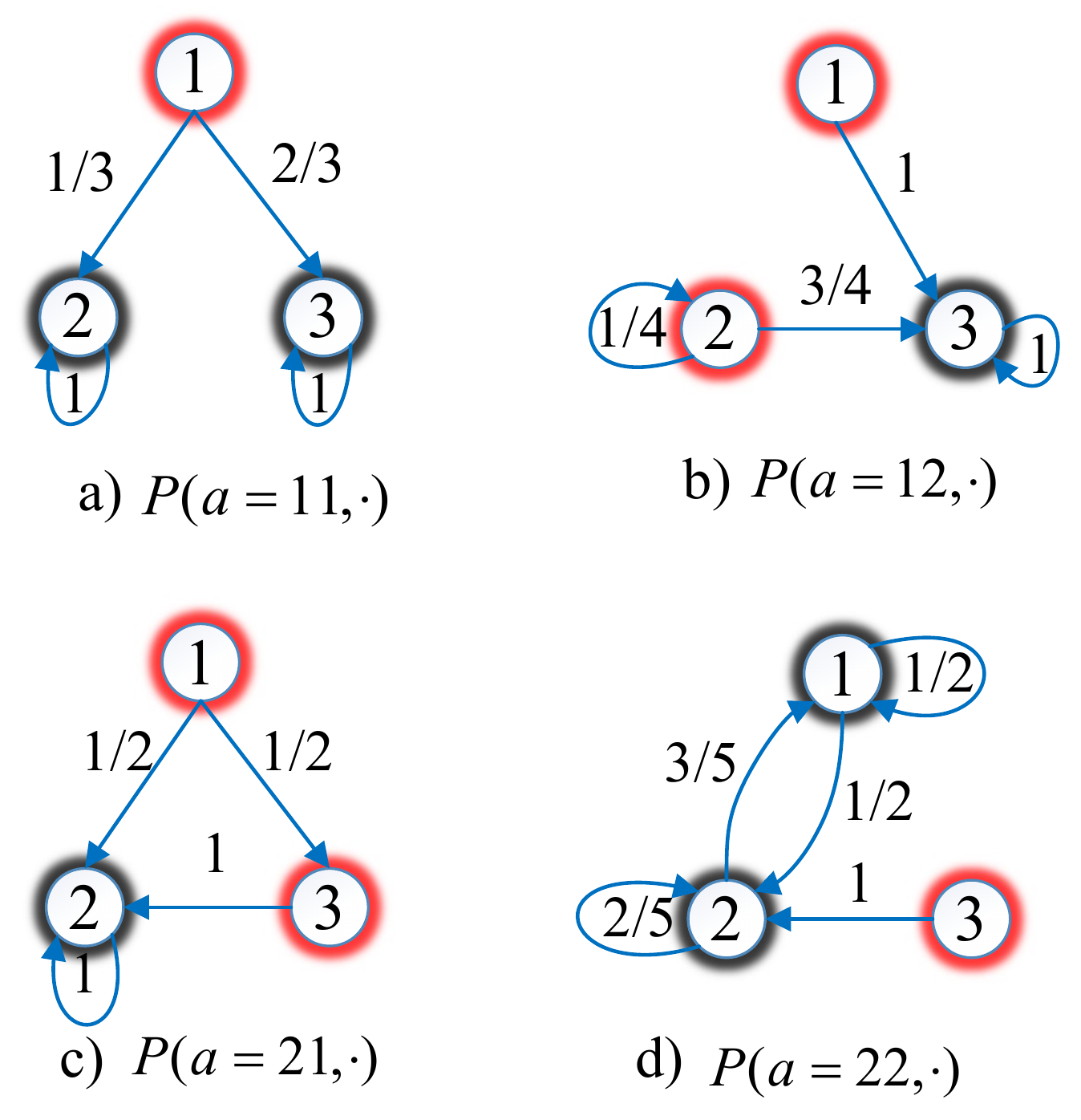}
    \caption{State Transition Diagram of Example \ref{ex1}}
    \label{fig1.1}  
\end{figure}

One can verify that the recurrent states of Markov chain $P(a=22,\cdot)$ is $x=1,x=2$, and $a=22$ is a pure Nash  equilibrium when $x=1,2.$ Therefore, action state pair $[a=22,x=1]$ and $[a=22,x=2]$ both are the recurrent state equilibria of Example \ref{ex1}, and $[a=22,x=1]\sim[a=22,x=2]$. Although $a=11$ is the pure Nash  equilibrium of $\mathcal{G}_1$, $x=1$ is a  transient state of Markov chain $P(a=11,\cdot)$. So $[a=11,x=1]$ is not a  recurrent state equilibrium.
\end{exa}
\subsection{State-based potential games}
State-based potential game, which is introduced by J. R. Marden \cite{jm12}, can guarantee the existence of a recurrent state equilibrium.
\begin{dfn}\label{dfn4} \emph{(State-based potential games)}\cite{jm12}
A state-based game $\mathcal{G}=\left\{N,\{A_i\},\{c_i\},X,P\right\}$ is called a state-based potential game if there exists a  function $\phi:A\times X\rightarrow \mathbb{R}$
such that for each  action state pair $[a,x]\in A\times X$, the following two conditions are satisfied:
\begin{enumerate}
\item For any agent $i\in N$ and action $a'_i\in A_i$
$$c_i(a'_i,a_{-i},x)-c_i(a,x)=\phi(a'_i,a_{-i},x)-\phi(a,x).$$

\item For any state $x'$ in the support of $P(a,x)$,
$$\phi(a,x')\geq \phi(a,x).$$
\end{enumerate}
$\phi$ is called a potential function of the state-based potential game, where $a_{-i}\in A_{-i}:=\prod_{j\neq i}A_j$ is the joint action profile other than agent $i$.
\end{dfn}

The first condition means that every state invariant game $G_x=\{N,A_i,c_i(\cdot,x)\}$ is a potential game.
The second condition ensures that any action state pair which maximizes the potential function is a recurrent state equilibrium of the state-based potential game. Denote by  $[a^*,x^*]$ the action state pair which maximizes the potential function, i.e., $[a^*,x^*]\in\arg\max_{[a,x]\in A\times X}\phi(a,x).$ Let $R(a|x)$ be the recurrent states of the Markov chain  $P(a,\cdot)$ starting from state $x$, which is by definition nonempty. Obviously, $R(a|x)\subseteq X(a|x).$ Therefore the second condition of Definition \ref{dfn4} can be relaxed as:\\
\emph{(2) If $[a^*,x^*]\in\arg\max_{[a,x]\in A\times X}\phi(a,x)$, then $[a^*,y]\in\arg\max_{[a,x]\in A\times X}\phi(a,x)$ for every $y\in R(a^*|x^*)$}.

\subsection{Learning in state-based games}
Roughly speaking, learning  in games is a decision-making process using available information. The difference of learning algorithm  between state-based games and normal form games is that for the former there is an additional  factor, \emph{state}, needed to be considered.

Consider a repeated  state-based game. The observed sequence  of agent $i$ at time $t$ is $\{\{a(\tau),x(\tau)\}_{\tau=0,1,\ldots,t-1},x(t)\}$. Let $O_i(t)$ denote the  obtained/available information of agent $i$ at time $t$, that is,
\begin{align*}
O_i(t):=\big\{~\{a(\tau),x(\tau)\}_{\tau=0,1,\ldots,t-1},x(t)\big\}.
\end{align*}

Generally speaking, the action updating mechanism of agent $i$ can be described by a response algorithm $f_i$ \cite{js93},
$$f_i:O_i(t)\rightarrow \Delta(A_{i}),$$
where $f_i$ is a function which maps  agent $i$'s available information $O_i(t)$  to a probability distribution over $i$'s own actions $A_i$.  Agent $i$ selects the action $a(t+1)\in A_i$ according to the probability distribution at time $t+1$. $\Delta(A_{i})$ denotes the set of probability distributions over $A_{i}$.

According to the available information used in making decisions, the most common learning algorithms can be categorized as \emph{uncoupled} learning algorithms and \emph{completely uncoupled} learning algorithms, whose definition are shown as follows.
\begin{dfn}\cite{ms13} A learning algorithm is called

i)  \emph{uncoupled} if the available information of agent $i$ used for decision-making is the payoff structure of himself and history sequence of the play, i.e.,
\begin{align*}
O_i(t)=\big\{~\{a(\tau),x(\tau)\}_{\tau=0,1,\ldots,t-1},x(t);~c_i(a,x)\big\}.
\end{align*}

ii) \emph{completely uncoupled} if the available information of agent $i$ used for decision-making is his own past realized payoffs and  actions, i.e.,
\begin{align*}
O_i(t)=\big\{~\{a_i(\tau),x(\tau),c_i(a(\tau),x(\tau))\}_{\tau=0,1,\ldots,t-1},x(t)\big\}.
\end{align*}

\end{dfn}
Replicator dynamics \cite{tb97}, best-reply \cite{pa18}, and fictitious play \cite{js05} are uncoupled learning algorithms. Regret learning \cite{sh00} and trial-and-error learning \cite{hp09} are completely uncoupled learning algorithms.

%

The  paper  focuses  on designing a \emph{natural} and \emph{effective} strategic learning algorithm  which converges to  recurrent state equilibrium of the state-based games. By \emph{natural} we require the  algorithm being uncoupled or completely uncoupled. By \emph{effective} we mean that the designed algorithm should converge to the equilibrium heuristically, not be trapped in an adjustment cycle, and not be predicted easily by each agent's opponents.
\section{A two-memory better reply learning algorithm}
\subsection{Available information}
Consider a repeated state-based game. Each agent seeks to maximize its myopic payoff. Agent  $i$ knows his own payoff function, but he doesn't know his opponents' ones. He can observe current state $x$ and his opponents' actions $a_{-i}\in A_{-i}$, but the agent doesn't know the structure of the Markovian state transition function $P$. Each agent can  recall the past $2$-period information, i.e. $2$-memory, at each time. Denote by $\xi_i(t)$ the information  used to make decision for agent $i$ at time $t\geq2$
$$\xi_i(t):=\big\{a(t-2),a(t-1),x(t);c_i(a,x)\big\}.$$
Then the response algorithm $f_i$ of agent $i$ has the following form
$$p_i(t)=f_i\big(\xi_i(t)\big)\in \Delta(A_i).$$

For any action state pair $[a,x]\in A\times X$, agent $i$'s \emph{strict better reply set} is defined as
$$B_i(a;x):=\big\{a'_i\in A_i:~c_i(a'_i,a_{-i},x)>c_i(a,x)\big\}.$$
For simplicity, let $B_i(t):=B_i(a(t-1);x(t)),~\forall t\geq 1$.

\subsection{The flow of the two-memory better reply learning algorithm}

Suppose the information of the past two periods at time $t\geq2$ is $[a(t-2),x(t-1)]\times [a(t-1),x(t)]\in (A\times X)\times(A\times X)$. The response algorithm $f_i$ of agent $i$ is defined as follows:
\begin{description}
   \item[(i)] Check whether $a(t-2)=a(t-1)$ or not at time $t$.
  \item[(ii)]  If $a(t-2)=a(t-1)$. Then each agent calculates $B_i(t)$ and check whether  $B_i(t)=\emptyset$ or not. If $B_i(t)=\emptyset$, then agent $i$ plays $a_i(t-1)$ next moment. Otherwise agent $i$ selects actions  according to a probability distribution on $A_i$, the support of which is $\{a_i(t-1)\}\cup B_i(t)$. Particularly, agent $i$ selects $a_i(t-1)$ with probability $\epsilon_i \in (0,1)$, the inertia of agent $i$, and the actions in $B_i(t)$ with equal probability.
  \item[(iii)] If $a(t-2)\neq a(t-1)$, then all agents take actions simultaneously according to their probability distributions with full support. Particularly, agent $i$ selects $a_i(t-1)$ with probability $\epsilon_i \in (0,1)$, and other actions in $B_i(t)$ with equal probability.
\end{description}

\begin{figure}[!hbtp]
\centering
\includegraphics[ width = 3.3 in]{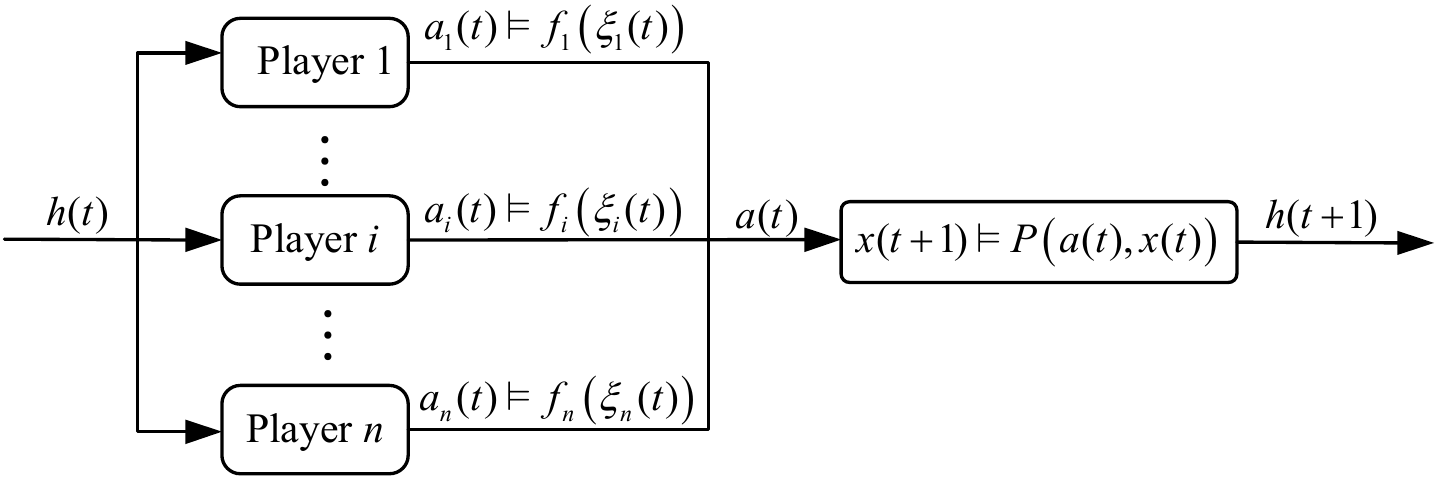}
\caption{Dynamics of State-based Games}
\label{fig3.1}
\end{figure}

Denote by $p_i^{a_i}(t)$  the probability that agent $i$ selects $a_i\in A_i$ at time $t$. The detailed algorithm of the proposed learning algorithm is shown in \textbf{Algorithm $1$}.
\begin{algorithm}[htbp]\label{alg1}
\caption{\textbf{.} Two memory better reply learning algorithm} 
\hspace*{0.02in} {\bf Input:} 
 $n$, $A_i$,  $c_i(a,x)$,  $X$,  $P(x,a)$, $\epsilon_i$.\\
\hspace*{0.02in} {\bf Output:} 
Recurrent state equilibrium of $\mathcal{G}$.
\begin{algorithmic}[1]
\State {\bf Initialization:} Choose a initial state $x(1)\in X$  randomly. Set simulation time $T\geq 3$. 
\For{$i=1:n$} 
¡¡¡¡\State $p_i^{a_i}(1)=\frac{1}{|A_i|}, \forall a_i\in A_i;$
\EndFor
\State {\bf end for}
\State $x(2)\vDash P(a(1),x(1));$
\For{$i=1:n$} 
¡¡¡¡\State $p_i^{a_i}(2)=\frac{1}{|A_i|}, \forall a_i\in A_i;$
\EndFor
\State {\bf end for}
\State $x(3)\vDash P(a(2),x(2));$
\For{$t=3:T$} 
¡¡¡¡\If{$a(t-2)=a(t-1)$} 
     \For{$i=1:n$}
       \If{$B_i(t)=\emptyset$}
¡¡¡¡        \State $a_i(t)=a_i(t-1);$
       \Else
¡¡¡¡       \State $p_i^{a_i(t-1)}(t)=\epsilon_i;$
           \State $p_i^{a_i}(t)=\frac{1-\epsilon_i}{|B_i(t)|}, \forall a_i\in B_i(t);$
       \EndIf
       \State {\bf end if}
     \EndFor
     \State {\bf end for}
¡¡¡¡\Else
¡¡¡¡¡¡¡¡\For{$i=1:n$} 
¡¡¡¡    \State $p_i^{a_i(t-1)}(t)=\epsilon_i;$
        \State $p_i^{a_i}(t)=\frac{1-\epsilon_i}{|A_i|-1}, \forall a_i\in A_i\setminus \{a_i(t-1)\};$
        \EndFor
        \State {\bf end for}
¡¡¡¡\EndIf
    \State {\bf end if}
\State $x(t+1)\vDash P(a(t),x(t));$
\EndFor
\State {\bf end for}
\State \Return
\end{algorithmic}
\end{algorithm}
\begin{rem}
The proposed learning algorithm is a $2$-memory, stochastic learning algorithm with inertia $\epsilon_i$ for agent $i$. It is a combination of \emph{testing}, \emph{searching}, and \emph{lock-in}.   Since the learning algorithm is $2$-memory, and every agent can can observe the opponents' actions. So each agents can tell whether $a(t-2)= a(t-1)$ or not. This is  \emph{testing}. The searching process consists of \emph{local search} and \emph{global search}. If $a(t-2)\neq a(t-1)$, then all agents take actions simultaneously according to their probability distributions with full support. This is a \emph{global stochastic search}, both for agents and actions.  If $a(t-2)= a(t-1)$ and $B_i(t)\neq\emptyset$, then agent $i$ will take actions from  $B_i(t)$. This is a \emph{local random search}.  If $a(t-2)= a(t-1)$ and $[a(t-2),x(t-2)]$ is an RSE, all agents will repeat their actions forever, which is called \emph{lock-in}.
\end{rem}
Denote by $h(t):=\{a(t-2),a(t-1),x(t)\}$ the past two plays, $t>2$. Then $\xi_i(t)=\{h(t);c_i(a,x)\},i\in N.$ 
The flow of the two-memory better reply learning algorithm can be described as in Fig. \ref{fig3.1}.
\subsection{Convergence of the proposed learning algorithm}

Consider a state-based game $\mathcal{G}=\big\{N,\{A_i\},\{c_i\},X,P\big\}$.  Let
\[
\bar{P}(\cdot, \cdot) := \frac{1}{|A|} \sum_{a \in A} P(a; \cdot, \cdot),
\]
and we know that $\bar{P}(\cdot, \cdot)  \in \mathbb{R}^{|X| \times |X|}$ is row stochastic. Then a Markov chain is defined by $\bar{P}$ with $X$ as its state space. Suppose $\mathcal{G}$ has at least one RSE, and let
$$A^*=\{a\in A| \text{there exists a state $x$, s.t. $[a,x]$ is a RSE}\}.$$
For $a \in A^*$, denote \[X(a) := \{x \in X :  \exists x^*\in X(a|x), \text{ s.t. } [a, x^*] \text{ is an RSE}\}.\]
The set $X(a),\forall a\in A^*$ contains all states from which the algorithm can reach an RSE class of action $a$ with positive probability by only adopting the same action $a$. Let $X^*:=\bigcup_{a \in A^*}X(a)\subseteq X$.

\begin{thm} \label{th3.7}
Consider a state-based game $\mathcal{G}=\left\{N,\{A_i\},\right.$ $\left.\{c_i\},X,P\right\}$, where the recurrent state equilibria exist. Suppose that either $X\setminus X^*  = \emptyset$, or $X\setminus X^* \not = \emptyset$ and the following assumptions hold:
\begin{description}
  \item[(i)]  For every recurrent class $\bar{R}$ of $\bar{P}$, there exists an action $a^* \in A$ and a state $x^* \in \bar{R}$ such that $[a^*, x^*]$ is an RSE.
  \item[(ii)]  $P(a; x, x) > 0$ for all $a \in A$ and $x \in X\setminus X^*$.
\end{description}
Then for any initial state $x_0 \in X$, if all agents play the game $\mathcal{G}$ by the proposed two memory better reply  learning algorithm,  the action state pair converges almost surely to an action invariant set of recurrent state equilibria.
\end{thm}
Conditions (i) and (ii) of Theorem \ref{th3.7} guarantee that there exists a positive probability ``path" which leads any initial action state pair to  an RSE. The proof of Theorem \ref{th3.7} is presented in the Appendix.

The following example shows that the assumption  (ii) of Theorems 7 avoids the situation where some desired actions cannot be selected according to the learning algorithm.

\begin{exa}\label{ex81}

Consider the following state-based game with $N = \{1, 2\}$, $A_1 = A_2 = \{\textup{C}, \textup{D}\}$, $X=\{1, 2, 3, 4\}$, and $A=\{CC,CD,DC,DD\}$. The payoff  bi-matrices  are shown in Table \ref{Tab3.4}-Table \ref{Tab3.7}.

\begin{table}[!htbp] 
\centering
\caption{Payoff Bi-Matrix for $x=1$ of Example \ref{ex81}   \label{Tab3.4}}
\begin{tabular}{ccc}
\hline Agent $1\backslash$Agent $2$ & \textup{C} & \textup{D}\\
\hline $\textup{C}$&$(5,4)$&$(2,3)$\\
 $\textup{D}$&$(4,2)$&$(3,1)$\\
\hline
\end{tabular}
\end{table}
\begin{table}[!htbp]
\centering
\caption{Payoff Bi-Matrix for $x=2$ of Example \ref{ex81}   \label{Tab3.5}}
\begin{tabular}{ccc}
\hline Agent $1\backslash$Agent $2$ &$\textup{C}$&$\textup{D}$\\
\hline $\textup{C}$&$(1,2)$&$(3,1)$\\
 $\textup{D}$&$(2,0)$&$(2,1)$\\
\hline
\end{tabular}
\end{table}
\begin{table}[!htbp]
\centering
\caption{Payoff Bi-Matrix for $x=3$ of Example \ref{ex81}  \label{Tab3.6}}
\begin{tabular}{ccc}
\hline Agent $1\backslash$Agent $2$ &$\textup{C}$&$\textup{D}$\\
\hline $\textup{C}$&$(-1,1)$&$(1,-1)$\\
 $\textup{D}$&$(1,-1)$&$(-1,1)$\\
\hline
\end{tabular}
\end{table}
\begin{table}[!htbp]
\centering
\caption{Payoff Bi-Matrix for $x=4$ of Example \ref{ex81}  \label{Tab3.7}}
\begin{tabular}{ccc}
\hline Agent $1\backslash$Agent $2$ &$\textup{C}$&$\textup{D}$\\
\hline $\textup{C}$&$(2,2)$&$(2,3)$\\
 $\textup{D}$&$(0,3)$&$(3,1)$\\
\hline
\end{tabular}
\end{table}

The  Markovian state transition matrices are as follows:
\[
P(\textup{CC}; \cdot, \cdot) = \begin{bmatrix}
1 & 0 & 0 & 0\\
0 & 1 & 0 & 0\\
0 & \frac12 & \frac12 & 0\\
0 & 0 & \frac12 & \frac12\\
\end{bmatrix}, \quad
P(\textup{CD}; \cdot, \cdot) = \begin{bmatrix}
1 & 0 & 0 & 0\\
\frac12 & \frac12 & 0 & 0\\
0 & 0 & 0 & 1\\
0 & 0 & 0 & 1\\
\end{bmatrix},
\]
\[
P(\textup{DC}; \cdot, \cdot) = \begin{bmatrix}
\frac12 & \frac12 & 0 & 0\\
0 & 0 & 0 & 1\\
0 & 0 & 0 & 1\\
0 & 0 & 0 & 1\\
\end{bmatrix}, \quad
P(\textup{DD}; \cdot, \cdot) = \begin{bmatrix}
1 & 0 & 0 & 0\\
0 & \frac12 & 0 & \frac12\\
0 & 0 & 0 & 1\\
0 & 0 & 0 & 1\\
\end{bmatrix}.
\]

It can be observed that the only RSE is $(CC, 1)$. Suppose that $x(0) = 4$, and the only possible choice of actions such that the system leaves the state $4$ and reaches the state $2$ is adopting \textup{CC} twice. This is because $a(0)$ must be $CC$ and $x(1) = 3$ with probability $1/2$. Although $a(1)$ can be any action in $A$, actions $\textup{CD}$, $\textup{DC}$, and $\textup{DD}$ make the system return to the state $4$. Therefore, $a(1)$ should be $\textup{CC}$ too, and $x(2) = 2$ with probability $1/2$ on the condition that $x(1) = 3$.

However, since $B_1(\textup{CC}, 2) = \{\textup{D}\}$ and $B_2(\textup{CC}, 2) = \emptyset$, the algorithm can only select actions from set $\{\textup{CC}, \textup{DC}\}$ at time $t = 2$. The choice $\textup{CC}$ makes the state of the system stay at $2$, while the latter makes $x(3) = 4$, and everything returns to the beginning. Thus, the algorithm cannot reach the RSE from the initial state $x(0) = 4$, though $\bar{P}$ is irreducible, and the assumptions (i) of Theorem \ref{th3.7} holds.
\end{exa}
\section{Existence of universal time-efficient learning algorithm}

\subsection{Time efficiency}
One may be interested in the complexity of the proposed learning algorithm, especially the time efficiency. The time efficiency of a learning algorithm is defined as follows:
\begin{dfn} \cite{ms13}
A learning algorithm is called \emph{time efficient} if the time for the algorithm to converge to an equilibrium is polynominal with respect to the number of agents.
\end{dfn}
\cite{sh10} proved that there does not exist any time-efficient uncoupled learning algorithm that converges to a pure Nash equilibrium for generic normal form games where such an equilibrium exists. As state-based games contain normal form games as its special case,  we can conclude that:

\begin{prp}
There does not exist any time-efficient uncoupled learning algorithms that converge to a recurrent state  equilibrium for general state-based games where such an equilibrium exists.
\end{prp}

\subsection{A counter example}
In fact, when it comes to state-based games, things  become  a bit more complicated. There is even no universal learning algorithm converging to a recurrent state equilibrium. We present the following example.
\begin{exa}\label{ex9} Consider the following state-based game with $N=\{1,2\},$ $A_1=A_2=\{1,2\},$ $X=\{1,2,3,4\}$. The payoff matrices are shown in Table \ref{Tab4.1}-Table \ref{Tab4.4}.
\begin{table}[!htbp] 
\centering
\caption{Payoff Bi-Matrix for $x=1$ of Example \ref{ex9}   \label{Tab4.1}}
\begin{tabular}{ccc}
\hline Agent $1\backslash$Agent $2$ &$1$&$2$\\
\hline $1$&$(5,~4)$&$(2,~3)$\\
 $2$&$(4,~2)$&$(3,~1)$\\
\hline
\end{tabular}
\end{table}
\begin{table}[!htbp]
\centering
\caption{Payoff Bi-Matrix for $x=2$ of Example \ref{ex9}   \label{Tab4.2}}
\begin{tabular}{ccc}
\hline Agent $1\backslash$Agent $2$ &$1$&$2$\\
\hline $1$&$(2,~2)$&$(3,~1)$\\
 $2$&$(0,~3)$&$(2,~1)$\\
\hline
\end{tabular}
\end{table}
\begin{table}[!htbp]
\centering
\caption{Payoff Bi-Matrix for $x=3$ of Example \ref{ex9}  \label{Tab4.3}}
\begin{tabular}{ccc}
\hline Agent $1\backslash$Agent $2$ &$1$&$2$\\
\hline $1$&$(-1,~1)$&$(1,~-1)$\\
 $2$&$(1,~-1)$&$(-1,~1)$\\
\hline
\end{tabular}
\end{table}
\begin{table}[!htbp]
\centering
\caption{Payoff Bi-Matrix for $x=4$ of Example \ref{ex9}  \label{Tab4.4}}
\begin{tabular}{ccc}
\hline Agent $1\backslash$Agent $2$ &$1$&$2$\\
\hline $1$&$(2,~2)$&$(2,~3)$\\
 $2$&$(0,~3)$&$(3,~1)$\\
\hline
\end{tabular}
\end{table}

The Markov transition matrices under different actions have the following form:
\begin{align*}
P(a;\cdot,\cdot)=
\begin{bmatrix}
p_{11}(a),&p_{12}(a),&0,&0\\
p_{21}(a),&p_{22}(a),&0,&0\\
0,&0,&p_{33}(a),&p_{34}(a)\\
0,&0,&p_{43}(a),&p_{44}(a)\\
\end{bmatrix},
\end{align*}
where $0<p_{ij}(a)<1$ is the probability that state $i$ transfers to state $j$, $\forall a\in\{11,12,21,22\}.$

It is obvious that action state pair $[a=11,x=1]$ and $[a=11,x=2]$ are  RSEs. For any learning algorithms, once the process enters action state pair $[a,x=3]$ or $[a,x=4]$, it  cannot escape from such an action state pair. Therefore, there does not exist any learning algorithms that converge to  a recurrent state equilibrium in such state-based games.
\end{exa}

According to Example \ref{ex9},   the following claim is obvious.

\begin{prp}
If for all Markov chain $P(a;\cdot,\cdot),\forall a\in A$, there exists a common closed set, denoted by $X^c\subseteq X$, s.t., such that, for all $x \in X^c$ and $a \in A$, $[a, x]$ is not an RSE. Then  there does not exist any uncoupled learning algorithm that converge to an RSE for generic state-based games even if such an equilibrium exists.
\end{prp}

The reason why there does not exist such learning algorithms  is that for a given state-based game the dynamic of the state $P(a;\cdot,\cdot)$ is pre-given, which is uncontrollable.
\section{Conclusion}

An extended model in game theory, called state-based games, is investigated in this paper. An uncoupled two memory learning algorithm is proposed. We proved that under certain reasonable conditions the proposed learning algorithm converges to a recurrent state equilibrium of a state-based games.  Since an additional degree of freedom is provided to help coordinate group behavior, state-based game is an  useful extended model in game-theoretic control. The existence of time-efficient universal learning algorithm is also investigated. A numerical example is presented to show that there is even no universal learning algorithm converging to a recurrent state equilibrium.
Future works will focus on the applications of the state-based game model and the learning algorithm to engineering control problems.

\appendix
  \renewcommand{\appendixname}{Appendix~\Alph{section}}

\textbf{Appendix}
\section{The  proposed algorithm and  corresponding Markov chain}
The proposed $2$-memory learning algorithm defines a discrete-time Markov chain $\{\omega(t), t \ge 0\}$ with finite state space $\Omega := X\times A \times X\times A\times X$, where $\omega(t) = [x(t),a(t),x(t + 1),a(t + 1),x(t + 2)]^T$, $t \ge 0$.

Let $x^{i} \in X$ and $a^{i} \in A$ be the state and action at time $i$, respectively. The initial distribution of the Markov chain $\{\omega(t)\}$ is
\begin{align*}
\begin{array}{cll}
\mathbf{Pr}\left\{\omega(0) = [x^{0},a^{0},x^{1},a^{1},x^{2}]^T\right\} \\
= \left(\prod_{1\le i \le n} \frac1{|A_i|}\right)^2 p(x^0) \text{P}(a^{0}; x^{0}, x^{1}) \text{P}(a^{1}; x^{1}, x^{2}),
\end{array}
\end{align*}
where $p : X \to [0,1]$ is the probability distribution of for initial state. For the sake of simplification, suppose the inertia of agent $i$ is the same, i.e., $\epsilon=\epsilon_i.$

Consider  any two states $\omega_1,\omega_2\in \Omega$ of the Markov chain $\{\omega(t)\}$, where $\omega_1 = [x^{1},a^{1},x^{2},a^{2},x^{3}]^T$ and $\omega_2 = [y^{1},b^{1},y^{2},b^{2},y^{3}]^T$. According to the learning algorithm, the transition probability from $\omega_1$ to $\omega_2$ of the Markov chain $\{\omega(t)\}$ is as follows:

\begin{enumerate}
  \item If $[y^{1},b^{1},y^{2}]\not= [x^{2},a^{2},x^{3}]$, then
        $$\mathbf{Pr}\left\{\omega(t+1)=\omega_2|\omega(t)=\omega_1\right\} =0.$$
  \item  If $[y^{1},b^{1},y^{2}] = [x^{2},a^{2},x^{3}]$ and $a^1\not =a^2$, then
        \begin{align*}
&~~\mathbf{Pr}\left\{\omega(t+1)=\omega_2|\omega(t)=\omega_1\right\} \\ &=\epsilon^{n-|H({b^1},{b^2})|}\cdot\prod\limits_{i\in H}\frac{1-\epsilon}{|A_{i}|-1}\cdot\textup{P}(b^{2}; y^{2}, y^{3}),
\end{align*}
where $H(a,b):=\{i\in N:a_i\neq b_i\},a,b\in A.$
  \item If $[y^{1},b^{1},y^{2}] = [x^{2},a^{2},x^{3}]$ and $a^1 =a^2$, then
        \begin{align*}
&~~\mathbf{Pr}\left\{\omega(t+1)=\omega_2|\omega(t)=\omega_1\right\} \\ &=\epsilon^{n-|H({b^1},{b^2})|-|N(b^1,y^2)|}\times\textup{P}(b^{2}; y^{2}, y^{3})\\
&~~\times\prod\limits_{i\in H}\frac{1-\epsilon}{|B_i(b^1,y^2)|}I_{B_i(b^1,y^2)}((b^2)_i),
\end{align*}
where $N(a,x):=\{i\in N: B_i(a,x)=\emptyset\},$ and $I_{B_i(a,x)}(b_i)$ is an indicator function such that $I_{B_i(a,x)}(b_i)= 1$ if $b_i \in B_i(a,x)$ and $I_{B_i(a,x)}(b_i) = 0$ if $b_i \notin B_i(a,x)$, $a\in A, x\in X, b_i\in A_i$.
\end{enumerate}

\section{Some lemmas used in proof of Theorem 7}
Denote $D(a, x) := \{b \in A : b_i \in B_i(a, x) \cup \{a_i\}, i\in N\}$ as the collection of action vectors whose entries are strict better reply actions for $a$ and $x$ or entries of $a$. From the definition, we know that $\{a\} \subseteq D(a, x)\subseteq A$ for any $a \in A$ and $x \in X$.

\begin{lem}\label{lem15}Consider a state-based game, where the RSE exists.
For any fixed initial value $x(0) = x^0$ and fixed action-state pairs $(a^0, x^1)$, $(a^1, x^2)$ of the learning algorithm, if there exists a positive integer $K \ge 2$ and a sequence of action-state pairs $\{(a^i, x^{i + 1}), ~2 \le i \le K\}$, where $a^i \in A, ~ x^{i+1} \in X, ~ 2 \le i \le K$, such that
\begin{description}
  \item[(i)] $P(a^2; x^2, x^3)P(a^3; x^3, x^4)\cdots P(a^{K}; x^{K}, x^{K + 1}) > 0$;
  \item[(ii)]  if $a^{k - 1} = a^{k}$ for some integer $k \in [1, K)$, then $a^{k + 1} \in D(a^{k}, x^{k + 1})$;
  \item[(iii)] $(a^{K}, x^{K + 1})$ is an RSE,
\end{description}
then the algorithm converges to some RSE  almost surely, by which we mean that $\textup{P}\{\tau < \infty\} = 1$, where $\tau := \min\{t \ge 2 : (a^t, x^{(t+1)}) \text{ is an RSE}~\}$, and, at the same time, that $a^{(\tau + t)} = a^{\tau}$, $x^{(\tau + t)} \in X(a^\tau|x^{(\tau+1)})$ for $t \ge 1$.
\end{lem}

\noindent{\it Proof}:
For convenience, let
$$\omega(t) :=[x^{t},a^{t},x^{t+1},a^{t+1},x^{t+2}]^T, \forall t\geq0,
$$
unless elsewhere stated.
The assumptions imply that, for any fixed initial state $\omega(0) =[x^{0},a^{0},x^1,a^1,x^2]^T$,
\[
\mathbf{Pr}\{\omega(K-1)|\omega(0)\} > 0.
\]
From the transition probability of $\{\omega(t)\}$ and that $(a^{K}, x^{K + 1})$ is an RSE, it follows that
\begin{align*}
&\mathbf{Pr}\{\omega(K+1) = [x^{K+1},a^{K},x^{K+2},a^{K},x^{K+3}]^T  \\
&~\qquad |\omega(K-1) = [x^{K-1},a^{K-1},x^K,a^K,x^{K+1}]^T\} > 0,
\end{align*}
where $x^{K+2}, x^{K+3} \in X(a^K|x^{K+1})$.

Thus,
\begin{align*}
&\mathbf{Pr}\{\omega(K+1) = [x^{K+1},a^{K},x^{K+2},a^{K},x^{K+3}]^T~~~\\
&~~~~~~~~~\qquad\qquad\quad |\omega(0) =[x^{0},a^{0},x^1,a^1,x^2]^T\} > 0,
\end{align*}
Therefore,  the algorithm can reach an RSE from any state  $\omega(0)\in\Omega$ with positive probability.
\hfill $\Box$

\begin{lem}\label{lem16}
Suppose that the following assumptions hold:
\begin{description}
  \item[(i)] $\bar{P}$ is irreducible;
  \item[(ii)] there exists an action $a^* \in A$ and a state $x^* \in X$ such that $(a^*, x^*)$ is an RSE;
  \item[(iii)] $P(a; x, x) > 0$ for all $a \in A$ and $x\in X$.
\end{description}
Then for any initial state $x \in X$, the algorithm converges to some RSE class a.s.
\end{lem}

\noindent{\it Proof}:
It suffices to validate the conditions in Lemma \ref{lem15} hold.

(i) For any fixed initial state $[x^{0},a^{0},x^1,a^1,x^2]$, if $a^0 \not= a^1$, and $(a^1, x^2)$ is an RSE, then the desired sequence of action-state pairs is obtained when we let $a^2 = a^1$. If $x^2 \in X(a^*|x^*)$, then let $a^2 = a^*$, and the desired sequence is obtained too.

Now assume that $a^0 \not= a^1$, that $(a^1, x^2)$ is not an RSE, and that $x^2 \not\in X(a^*|x^*)$. From assumption (i), it follows that, for $x^2 \in X$, there exists a positive integer $K_1 \ge 3$ such that
\[
\bar{P}(x^2, x^3)\bar{P}(x^3, x^4)\cdots \bar{P}(x^{K_1 - 1}, x^{K_1}) > 0,
\]
where $x^i \not= x^*$, $2 \le i < K_1$, and $x^{K_1} = x^*$.
The definition of $\bar{P}$ implies that there exists a sequence of action-state pairs $\{(a^{i}, x^{i + 1}), ~2 \le i < K_1\}$ such that
\[
P(a^2;x^2, x^3)P(a^3;x^3, x^4)\cdots P(a^{K_1 - 1}; x^{K_1 - 1}, x^*) > 0,
\]
where $x^i \not= x^*$, $2 \le i < K_1$. Let $a^{K_1} = a^*$.

Without loss of generality, suppose that $(a^i, x^{i + 1})$ is not an RSE for all $2 \le i < K_1$. Otherwise let $
\tilde{K}_1 := \min\{2 \le i < K_1 : (a^i, x^{i + 1}) \text{ is an RSE}\}$ and consider the sequence $\{(a^{i}, x^{i+1}), ~0 \le i \le \tilde{K}_1\}$.

Suppose that there exists some integer $k \in [1, K_1)$ such that $a^{k - 1} = a^{k}$ but $a^{k + 1} \not \in D(a^k, x^{k + 1})$. 
Denote $\hat{k} := 1 + \max\{t \in [0, k - 1) : a^t \not= a^{k - 1}\}$. The assumption $a^0 \not= a^1$ implies that $\hat{k} \ge 1$. Insert an action $\tilde{a}^{i} \not = a^{i}$ between $a^{i}$ and $a^{i + 1}$, $\hat{k} \le i < k$. In fact, $\tilde{a}^{i}$, $\hat{k} \le i < k$, can be the same action vector. Assumption (iii) ensures that
\begin{align*}
&P(a^{\hat{k}}; x^{\hat{k}}, x^{\hat{k} + 1}) P(\tilde{a}^{\hat{k}}; x^{\hat{k} + 1}, x^{\hat{k} + 1}) P(a^{\hat{k} + 1}; x^{\hat{k} + 1}, x^{\hat{k} + 2})\cdots\\
&P(a^{k - 1}; x^{k - 1}, x^{k}) P(\tilde{a}^{k - 1}; x^{k}, x^{k}) P(a^{k}; x^{k}, x^{k + 1}) > 0.
\end{align*}
The condition (ii) in Lemma \ref{lem15} is satisfied for this new sequence of action-state pairs, and the desired sequence  is obtained in this way.

(ii) If $a^0 = a^1$, and $(a^1,x^2)$ is an RSE, then let $a^2 = a^1$ and $x^3 \in X(a^1|x^2)$.

(iii) If $a^0 = a^1$, but $(a^1,x^2)$ is not an RSE, then, according to the learning rule, one can choose $a^2 \not=a^1$. By applying the argument above to $(x^1,a^1,x^2,a^2,x^3)$, we can obtain the desired sequence of action-state pairs.
\hfill $\Box$

\begin{lem}\label{lem17}
Suppose that the following assumptions hold:\\
(i) for every recurrent class $\bar{R}$ of $\bar{P}$, there exists an action $a^* \in A$ and a state $x^* \in \bar{R}$ such that $(a^*, x^*)$ is an RSE;\\
(ii) $P(a; x, x) > 0$ for all $a \in A$ and $x\in X$.\\
Then for any initial state $x \in X$, the algorithm converges to some RSE class a.s.
\end{lem}

\noindent{\it Proof}:
From the proof of Lemma \ref{lem16}, it suffices to show that the conditions in Lemma \ref{lem15} still hold when $a^0 \not= a^1$, and $x^2$ is a transient state of $\bar{P}$. If there exists an action $a^* \in A$ such that $(a^*, x^2)$ is an RSE, then let $a^2 = a^*$ and the desired sequence is obtained. Otherwise, since $x^2$ is transient for $\bar{P}$, we know that there exists a positive integer $K_1 \ge 3$ and a recurrent state of $\bar{P}$, $\tilde{x}$, such that
\[
\bar{P}(x^2, x^3)\bar{P}(x^3, x^4)\cdots \bar{P}(x^{K_1 - 1}, x^{K_1}) > 0,
\]
where $x^i \not= \tilde{x}$, $2 \le i < K_1$; $x^{K_1} = \tilde{x}$; $(\tilde{a}, \tilde{x})$ is an RSE for some $\tilde{a} \in A$. The definition of $\bar{P}$ implies that there exists a sequence of action-state pairs $\{(a^{i}, x^{i + 1}), ~2 \le i < K_1\}$ such that
\[
P(a^2;x^2, x^3)P(a^3;x^3, x^4)\cdots P(a^{K_1 - 1}; x^{K_1 - 1}, \tilde{x}) > 0,
\]
where $x^i \not= \tilde{x}$, $2 \le i < K_1$. Let $a^{K_1} = \tilde{a}$.

We can obtain the desired sequence by applying the same argument in Lemma \ref{lem16}.
\hfill $\Box$

\section{The proof of Theorem 7}

\noindent{\it Proof}: Before proving the theorem, we point out the following facts: if
\begin{enumerate}
  \item[(a)] the action state pair $[a(t),x(t)]$ is a RSE,
  \item[(b)] the action $a(t)$ is repeated for the next time, i.e. $a(t+1)=a(t)$,
  \item[(c)] all agents use the proposed learning algorithm,
\end{enumerate}
 then for $\forall t'> t+1$, $[a(t'),x(t')]$ will be a RSE and $a(t')=a(t).$ Therefore according to Lemma \ref{lem15}, the proof of Theorem \ref{th3.7} is  equivalent to proving the following statements: for any action state pair $[a(t),x(t)], \forall t>0,$ there exists a finite timesteps $T>0$ and a positive probability $\rho\in (0,1]$ such that $[a(t+T),x(t+T)]$ is a RSE and $a(t+T)=a(t+T+1)$ with at least probability $\rho>0$.


Denote  by $S:=(A\times X)\times(A\times X)$. Split  $S$ into four disjoint parts:
\begin{align*}
\begin{array}{cll}
S_1&:=&\big\{[a,x]\times[b,y]\in S: [a,x]\sim[b,y]\big\};\\
S_2&:=&\big\{[a,x]\times[b,y]\in S: [a,x]\nsim[b,y] \text{ and $[b,y]$ is a}\\
&~&~~\text{RSE}\big\};\\
S_3&:=&\big\{[a,x]\times[b,y]\in S: [a,x]\nsim[b,y], \text{$[b,y]$ is not a}\\
&~&~~\text{RSE, and $a\neq b$}\big\};\\
S_4&:=&\big\{[a,x]\times[b,y]\in S: [a,x]\nsim[b,y], \text{$[b,y]$ is not a}\\
&~&~~\text{RSE, and $a=b$}\big\}.
\end{array}
\end{align*}

Before starting the proof, we suppose $\epsilon_i=\epsilon,\forall i\in N.$ This assumption  will not affect the results.

\emph{Case $1$:} Suppose the play of the past two periods at time $t>2$ is $[a(t-2),x(t-1)]\times [a(t-1),x(t)]\in S_1$. Then there exists a recurrent state equilibrium set $R(a,x)$ such that $[a(t-2),x(t-1)]\times [a(t-1),x(t)]\in  R(a,x)\times R(a,x)$. It follows that $a(t-2)=a(t-1)=a$. According to the proposed learning algorithm, for any $t'\geq t$,  $[a(t')=a,x(t'+1)]$ is a RSE, and we are done.

\emph{Case $2$:} Suppose  $[a(t-2),x(t-1)]\times [a(t-1),x(t)]\in S_2$. Denote by
\begin{align*}
\begin{array}{cll}
S_2^1&:=&\big\{[a,x]\times[b,y]\in S_2:  a=b \big\},\\
S_2^2&:=&\big\{[a,x]\times[b,y]\in S_2:  a\neq b \big\}.\\
\end{array}
\end{align*}

\begin{itemize}
\item
If $[a(t-2),x(t-1)]\times [a(t-1),x(t)]\in S_2^1$, then according to condition (i) of the proposed learning algorithm, all agent will take $a(t-1)$ at time $t$ with probability $1$. So the action state pair $[a(t)=a(t-1),x(t+1)]$ is a RSE. Therefore $[a(t-1),x(t)]\times [a(t),x(t+1)]\in S_1.$ According to the above argument, we are done.
\item
If $[a(t-2),x(t-1)]\times [a(t-1),x(t)]\in S_2^2$, then according to condition (ii) of the proposed learning algorithm, all agent will take action simultaneously. The probability of $a(t)=a(t-1)$ is at least $\epsilon^n$. Hence the probability of $[a(t-1),x(t)]\times [a(t),x(t+1)]$ transfers into $S_1$ after $2$ steps with at least probability $\epsilon^n$.
\end{itemize}
Once $[a(t-1),x(t)]\times [a(t),x(t+1)]$ transfers into $S_1$, it will stay in a recurrent state equilibrium set $R(a,x)$ forever.

\emph{Case $3$:} Suppose  $[a(t-2),x(t-1)]\times [a(t-1),x(t)]\in S_3$. Let $[a^*,x^*]$ be an RSE of $\mathcal{G}$. According to algorithm (ii) of the proposed learning algorithm, all agents will take action simultaneously. The probability of $a(t)=a^*$ is
$$\delta_1=\epsilon^{n-|H(a(t-1),a^*)|}\prod\limits_{i_j\in H}\frac{1-\epsilon}{|A_{i_j}|-1}>0,$$
where $H(a(t-1),a^*)=\{i:a_i(t-1)\neq a_i^*\}$.
\begin{itemize}
  \item If $x(t+1)=x'\in X(a^*|x^*)$. Denote by $\gamma_1>0$ the probability that $x(t)\rightarrow x(t+1)=x'$ under the action $a^*.$ Then $[a(t-1),x(t)]\times [a(t),x(t+1)]$ transfers into $S_2$ with probability $\delta_1\cdot\gamma_1>0.$
  \item If $x(t+1)\notin X(a^*|x^*)$.
(i) If $X\setminus X^*=\emptyset.$ According to the definition of $X^*$, we know that there exists an action $a$ such that $[a,x(t)]$ is an RSE. As $a(t-2)\not= a(t-1)$, according to the learning algorithm the probability of $a(t)=a$ is positive. And $x(t+1)\in X(a(t)|x(t))$. The probability of $a(t)=a$ is
$$\delta_2=\epsilon^{n-|H(a(t-1),a)|}\prod\limits_{i_j\in H}\frac{1-\epsilon}{|A_{i_j}|-1}>0,$$
where $H(a(t-1),a)=\{i:a_i(t-1)\neq a_i\}$.
Then $[a(t-1),x(t)]\times [a(t)=a,x(t+1]$ transfers into $S_2$ with  probability $\delta_2\cdot P(a;x(t),x(t+1))>0$. (ii) If $X\setminus X^*\neq\emptyset$, and there exists an action $b\in A^*$, such that $x(t)\in X(b)$. As $a(t-2)\not= a(t-1)$, according to the learning algorithm, let $a(t)=b$. According to the definition, we know that there exists a finite integer $K>0$ and a state $x\in X$ such that
\begin{align}\label{eq1}
x(t)\stackrel{b}{\longrightarrow}~\stackrel{b}{\cdots}~\stackrel{b}{\longrightarrow}x(t+K)\in X(b|x),
\end{align}
where $(b,x)$ is an RSE. Then $[a(t+K-1)=b,x(t+K)]\times [a(t+K)=b,x(t+K+1)]$ transfers into $S_2$ with probability $\epsilon^{nK}\cdot\delta_3\cdot P(b;x(t+K),x(t+K+1))>0$, where
$$\delta_3=\epsilon^{n-|H(a(t+K-2),b)|}\prod\limits_{i_j\in H}\frac{1-\epsilon}{|A_{i_j}|-1}>0.$$
(iii) If $X\setminus X^*\neq\emptyset$, and  $x(t)\notin X^*$.
From the proof of Lemma \ref{lem16} and Lemma \ref{lem17}, we know that there exists a positive integer $K_1 \ge 0$ and a recurrent state $\tilde{x}$ of $\bar{P}$, such that
\[
\bar{P}(x(t), x(t+1))\cdots \bar{P}(x({t+K_1-1}), x({t+K_1})) > 0,
\]
where $x(t+\tau) \notin X^*$, $0 \leq\tau < K_1$, $x({t+K_1}) = \tilde{x}$ and $\tilde{x}\in X(\tilde{a})$ for some $\tilde{a}\in A^*$. Moreover, the definition of $\bar{P}$ implies that there exists a sequence of action-state pairs $\{(a(t+\tau), x(t+\tau+1)), ~0 \leq\tau < K_1\}$ such that
\begin{align*}
&P(a(t);x(t), x(t+1))P(a(t+1);x(t+1), x(t+2))\\
&\cdots P(a(t+K_1 - 1); x(t+K_1 - 1), \tilde{x}) > 0,
\end{align*}
where $x(t+\tau) \notin X^*$, $0 \leq\tau < K_1$, and $\tilde{x}\in X(\tilde{a})$ for some $\tilde{a}\in A^*$.
Assumption (ii) in Theorem $7$ ensures that by applying the same argument as in Lemma \ref{lem16} and Lemma \ref{lem17}, with a slight abusement of notations, we can obtain a sequence of action-state pairs $\{(a(t+\tau), x(t+\tau+1)), ~0 \leq\tau \leq K_1\}$,  which satisfy the all the conditions in Lemma \ref{lem15}. Using the same arguments in above condition (ii), there exists a  a finite integer $K_2>0$, such that
\begin{align*}
&[a(t+K_1+K_2-1),x(t+K_1+K_2)]\\
&~~~~~~~~~~\times [a(t+K_1+K_2),x(t+K_1+K_2+1)]
\end{align*}
transfers into $S_2$ with positive probability.

\end{itemize}
 According to the arguments in Case $2$, we can conclude that any state in $S_3$ will transfer into $S_1$ with a positive probability after finite steps.

\emph{Case $4$:} Suppose  $[a(t-2),x(t-1)]\times [a(t-1),x(t)]\in S_4$. Let $[a^*,x^*]$ be an RSE of $\mathcal{G}$.
\begin{itemize}
  \item If $a(t-1)=a^*,$ according to the arguments in Case 3, we can conclude that
  $$[a(t-1)=a^*,x(t)]\times [a(t)=a^*,x(t+1)]$$
   will transfer into $S_1$ will a positive probability after finite steps. Similar with the arguments in Case $2$, the probability that $[a(t-1),x(t)]\times [a(t),x(t+1)]$ transfers into $S_2$ is at least $\delta_2\cdot\gamma_3\cdot\epsilon^{nm}>0$.
   \item If $a(t-1)\neq a^*,$ there must be an agent $i\in N$ with an action $a_i'\in A_i$ for some state $x'\in X(a(t-1)|x(t))$ such that
   $$c_i(a_i',a_{-i}(t-1),x')>c_i(a(t-1),x').$$
   Otherwise, $a(t-1)=a^*.$ Since $x'\in X(a(t-1)|x(t)),$   there exists a time $t'\in \{t+1,\ldots,t+m+1\}$ such that
\begin{align}\label{eq2}
\textbf{Pr}[x(t')=x']>\theta>0
\end{align}
conditioned on the events $x(t)$, $a(t-1)=a(t)=\cdots=a(t'-1).$
The above events  happen with at least probability $\theta\cdot\epsilon^{nm}$. Denote by $a'=(a_i',a_{-i}(t-1))$. If $(a',x')$ is an RSE.
Then
$$[a(t'-1)=a(t-1),x(t)]\times [a',x']$$
transfers into $S_2.$
Notice that $a'\neq a(t'-1).$ If $(a',x')$ is not an RSE,
$$[a(t'-1)=a(t-1),x(t)]\times [a',x']$$
transfers into $S_3.$ According to the arguments in Case $2$ and Case $3$, we can conclude that any state in $S_4$ will transfer into $S_1$ will a positive probability after finite steps.
\end{itemize}
Summarizing  Case $1$, Case $2$, Case $3$ and Case $4$, we conclude that  for any fixed initial state $x(0)$ and any action state pair $[a(t),x(t+1)], \forall t>0,$ there exists a finite time $T$  such that $[a(t+T),x(t+T+1)]$ is an RSE and $a(t+T)=a(t+T+1)$ with a positive probability. 
%

\hfill $\Box$

\end{document}